\documentclass[preprint,1p,times]{elsarticle}

\usepackage{graphicx} 
\usepackage{xcolor}
\usepackage{ulem}
\usepackage{amsmath}
\usepackage{booktabs}
\usepackage{tabularx}
\usepackage{geometry}
\usepackage{algorithm}
\usepackage{algorithmic}
\usepackage{ragged2e}
\usepackage{xspace}
\usepackage{tikz}
\usetikzlibrary{shapes,snakes}
\usepackage{comment}

\newtheorem{example}{Example}[section]

\newcommand{\supplycenters}{supply warehouses\xspace}
\newcommand{\fc}{SW\xspace}

\newcommand{\sortcenters}{sort centers\xspace}
\newcommand{\scc}{SC\xspace}
\newcommand{\distributioncenter}{distribution station\xspace}
\newcommand{\distributioncenters}{distribution stations\xspace}
\newcommand{\ds}{DS\xspace}
\newcommand{\path}{path\xspace}
\newcommand{\paths}{paths\xspace}
\newcommand{\tdtime}{truck departure time\xspace}
\newcommand{\tdtimes}{truck departure times\xspace}
\newcommand{\tdt}{TDT\xspace}
\newcommand{\tdts}{TDTs\xspace}
\newcommand{\obj}{$F_{\text{black box}}$\xspace}

\begin{document}

\begin{frontmatter}
\title{Optimization of Next-Day Delivery Coverage using Constraint Programming and Random Key Optimizers}

\author[1]{Kyle Brubaker\corref{cor1}}
\ead{johbruba@amazon.com}
\cortext[cor1]{Corresponding author}
\author[1]{Kyle E. C. Booth}
\author[1]{Martin J. A. Schuetz}
\author[2]{Philipp Loick}
\author[2]{Jian Shen}
\author[2]{Arun Ramamurthy}
\author[2]{Georgios Paschos}
\address[1]{Amazon Advanced Solutions Lab, Seattle, Washington 98170, USA}
\address[2]{Amazon EU Science and Technology, Luxembourg}

\begin{abstract}
We consider the logistics network of an e-commerce retailer, specifically the so-called ``middle mile'' network, that routes inventory from \supplycenters to \distributioncenters to be ingested into the terminal (``last mile'') delivery network. 
The speed of packages through this middle mile network is a key determinant for the ultimate delivery speed to the end user. 
An important target for a retailer is to maximize the fraction of user orders that can be serviced within one day, i.e., next-day delivery. 
As such, we formulate the maximization of expected next-day delivery coverage within the middle-mile network as an optimization problem, involving a set of temporal and capacity-based constraints on the network and requiring the use of a black-box model to evaluate the objective function.
We design both exact constraint programming (CP) and heuristic random-key optimizer (RKO) approaches, the former of which uses a proxy objective function. 
We perform experiments on large-scale, real-world problem instances and show that both approaches have merit, in that they can match or outperform the baseline solution, a bespoke greedy solver with integrated local search, in expected next-day delivery coverage.
Our experiments focus on two high-level problem definitions, starting with a base problem and then adding more complexity, and also explore the generalization of the solvers across a range of problem instance sizes.
We find that a hybrid model using RKO and a bespoke local search protocol performs best on the full problem definition with respect to expected next-day delivery (increase of +50 basis points [bps] over baseline) but can take days to run, whereas the hybrid model using CP and local search is slightly less competitive (+20 bps) but takes only hours to run. 
\end{abstract}

\begin{keyword}
Middle mile logistics \sep Random key optimization \sep Constraint programming 


\end{keyword}

\end{frontmatter}

\section{Introduction}\label{sec:introduction}

For an e-commerce retailer, such as Amazon or Alibaba, it is important to maximize the fraction of customer orders that can be serviced by the next business day, i.e., next-day delivery.
It was estimated in the context of Alibaba that sales would increase by 13.3\% if three-day deliveries were moved towards two-day \citep{deshpande2023logistics}. 
Thus, retailers invest significant effort into optimizing their logistics networks for delivery speed. 
These delivery networks can be split into three distinct components: first mile, middle mile, and last mile. 
First mile represents the sourcing of products from suppliers to the retailer’s logistics network that generally happens prior to customer orders, meaning the first mile does not impact delivery speed of items to customers.
The last mile comprises shipments from \distributioncenters to the customer location, e.g., the customer’s home or office. 
It is generally assumed that if products can be injected into the last mile network prior to the appropriate cutoff, they will be delivered on the same day as their injection.
The middle mile then represents the intermediate network, connecting the \supplycenters, where the e-commerce retailers store the products, to the \distributioncenters located in proximity to end customers. 
The delivery speed from customer order to delivery is primarily driven by the \tdtimes in the middle-mile network.
Therefore, in this work, we seek to optimize the delivery speed within the middle mile network.

The middle mile network consists of origin nodes, intermediate nodes, and terminal nodes, with \paths (consisting of one or more legs) connecting origins to destinations either directly or through a set of intermediate nodes. 
The origin nodes (\supplycenters (\fc)) receive products from suppliers, and the destination nodes (\distributioncenters (\ds)) distribute the packages across appropriate local delivery vehicles to reach customers. 
The intermediate nodes (\sortcenters (\scc)) allow for the consolidation of volume throughout the network, increasing efficiency and reducing waste. 
In our work, we seek to optimize the so-called \tdtimes (\tdt) along each leg of this network, such that expected next-day delivery coverage is maximized. 

Finding good \tdt values is a difficult combinatorial optimization problem due to the size of the network, overlapping product coverage \citep{benidis:23}, numerous problem constraints, and the complex objective function. 
Specifically, \tdt values must satisfy node shift and labor constraints, capacity constraints, as well as node-specific time windows. 
Further, once a set of \tdt values have been determined, the quality of the solution must be evaluated using an expensive external black box function. 
Finally, for a given \path, demand coverage increases as the \tdt is moved later, but after some point there is a delivery day transition, i.e. no longer next-day, and so the solution reward must balance increasing demand coverage and minimizing delivery day.
This precludes a more straightforward modeling of the problem with paradigms such as integer linear programming which rely on analytically (e.g., linearly) expressible objective functions. 

A baseline approach for this problem involves a simple, yet powerful greedy search heuristic algorithm with a subsequent local search (i.e., the ``baseline'') as described in Section \ref{sec:methods-incumbent}.
In this work we explore whether alternative formulations and solvers can improve upon this baseline approach. 
We split our efforts into two distinct tracks: (i) using exact optimization solvers based on constraint programming (CP) \citep{rossi2006handbook} leveraging a proxy objective function, and (ii) using heuristic optimization solvers based on the random-key optimization (RKO) formalism \citep{chaves:24rko, chaves24, londe:24a, londe:24b, schuetz:22}. 

Constraint programming (CP) is an exact optimization method based on the model-and-solve paradigm where problems are first modeled with decision variables, constraints, and an objective function, and then constraint solvers (typically) employ a branch-and-infer search to find satisfying solutions to the modeled problem \citep{rossi2006handbook}. 
It is a popular alternative to integer programming and has been successfully used to solve a wide variety of problems, including those in scheduling \citep{baptiste2001constraint}, vehicle routing \citep{booth2019constraint}, and network design \citep{barra2007solving}. 
The CP paradigm emphasizes logical inference and is particularly effective when variable domains (the set of values a particular variable can be assigned) can be significantly pruned throughout the search. 
It is not straightforward, however, to incorporate the black box function calls during the CP search, and so we leverage a proxy objective function designed to approximate the more complex objective evaluation as described in Section \ref{sec:problem}. 

Random-key optimization (RKO) is a flexible search-based optimization framework that splits the optimization routine into two key components: the problem-independent search and the problem-dependent decoder. 
The problem-dependent decoder accepts a random key vector of some dimensionality $A$ ($\chi \in [0,1)^{A}$) and maps that to a candidate problem solution in the original problem space, calculates a fitness (or, equivalently, cost) for that solution, and provides that fitness signal $f(\chi)$ back to the problem-independent searcher to inform the next algorithmic step. 
The search itself can be performed by a variety of combinatorial algorithms (e.g., genetic algorithms, physics-inspired algorithms such as simulated annealing, or swarm-based methods such as particle-swarm optimization) \citep{chaves:24rko}, and may or may not require gradient information to take steps in this space. 
This framework is attractive as the black box objective can be easily called to update the algorithm throughout the search process.

\paragraph{Contributions}

Our work expands on \citep{benidis:23} in terms of both problem scope and algorithm development. 
First, we extend the network to include intermediate \scc nodes, increasing the number of variables to solve for and the number of constraints that must be considered.
Beyond this, we consider a number of additional or modified constraints. 
First, we add rolling time window capacity constraints, ensuring the number of trucks leaving a node, or arriving at a node, within the time window are limited. 
Second, we consider labor capacity in the \ds nodes, and ensure that the arriving volume is sufficient to avoid idle times and keep the labor fully engaged during the work shifts. 
Third, we expand the limit on the number of trucks, i.e., Constraint 3c in \citep{benidis:23}, allowing for multiple trucks on a subset of legs.
Finally, the application of heuristic RKO-based optimizers is novel to this problem, as is the application of CP to the problem. 
In particular, the CP approach uses a novel reformulation of the problem, as detailed in Section \ref{sec:methods-cp}.

\paragraph{Paper Structure}

The remainder of this paper is organized as follows.
In Section~\ref{sec:problem} we first define the problem, including the network, decision variables, objective, and constraints.
Then in Section~\ref{sec:methods} we describe our methodologies, including converting a solution candidate into a network-level plan, details of how the RKO and CP solvers work, the baseline solution, and the hybrid approaches considered.
In Section~\ref{sec:experiments} we present the results of these approaches as applied to the full network, considering the different problem versions described in Section~\ref{sec:problem}.
Finally, in Section~\ref{sec:conclusion} we present a discussion of the results, as well as learnings from the work.

\section{Problem definition}
\label{sec:problem}

\begin{figure}[h]
  \centering
  \begin{tikzpicture}[
    node distance=4cm,
    box/.style={rectangle,draw,minimum width=1cm,minimum height=1cm},
    diamondshape/.style={diamond,draw,minimum width=1cm,minimum height=1cm},
    circleshape/.style={circle,draw,minimum width=1cm}
]

\node[anchor=south] at (0,-1.15) {Supply};
\node[anchor=south] at (0,-1.6) {Warehouses (SW)};
\node[anchor=south] at (4,-1.25) {Sort Centers (SC)};
\node[anchor=south] at (8,-1.15) {Distribution};
\node[anchor=south] at (8, -1.6) {Stations (DS)};

\node[box] (SW1) at (0,3) {SW1};
\node[box] (SW2) at (0,1) {SW2};

\node[diamondshape] (SC1) at (4,3) {SC1};
\node[diamondshape] (SC2) at (4,1) {SC2};

\node[circleshape] (DS1) at (8,4) {DS1};
\node[circleshape] (DS2) at (8,2) {DS2};
\node[circleshape] (DS3) at (8,0) {DS3};

\node at (1.9,3.3) {$\delta_{(SW1,SC1)}$};

\draw[-stealth] (SW1) -- (SC1);
\draw[-stealth] (SW1) -- (SC2);
\draw[-stealth] (SW2) -- (SC1);
\draw[-stealth] (SW2) -- (SC2);
\draw[-stealth] (SW2) to[out=-20,in=-180]  (DS3);

\draw[-stealth] (SC1) -- (SC2);

\draw[-stealth] (SC1) -- (DS1);
\draw[-stealth] (SC1) -- (DS2);
\draw[-stealth] (SC2) -- (DS1);
\draw[-stealth] (SC2) -- (DS2);
\draw[-stealth] (SC2) -- (DS3);

\end{tikzpicture}
  \vspace*{-0.75cm}
  \caption{Example middle mile network with $|\mathcal{V}|=7$ and $|\mathcal{A}|=11$. The \supplycenters (\fc) serve as \path origins, connecting to downstream \sortcenters (\scc) or \distributioncenters (\ds). The travel time between nodes $i$ and $j$ is denoted as $\delta_{(i,j)}$. The network has both direct \paths (e.g., $(SW2,DS3)$) and \paths with multiple legs (e.g., $((SW1, SC1), (SC1, DS1))$.}
  \label{fig:problem_schematic}
\end{figure}
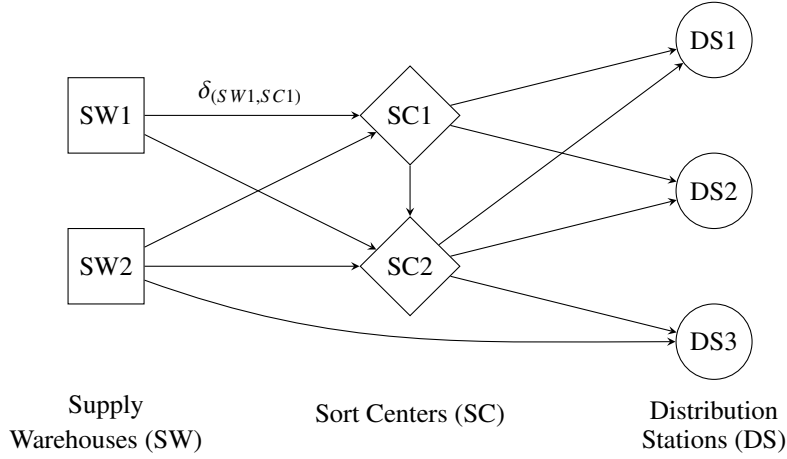

Let the middle mile network be defined on a directed graph $\mathcal{G} = (\mathcal{V}, \mathcal{A})$ where $\mathcal{V}$ is the set of $n$ nodes (locations) in the network. Further, let $\mathcal{F}$, $\mathcal{D}$, $\mathcal{S} \subset \mathcal{V}$ be the sets of \supplycenters (\fc), \distributioncenters (\ds), and \sortcenters (\scc), respectively.  Finally, $\mathcal{A}$ is the set of $m$ arcs (legs) in the network, and each element $(i,j) \in \mathcal{A}$ represents a leg with origin $i$ and destination $j$. A small example middle mile network is depicted in Figure \ref{fig:problem_schematic}. 

The transit time for a leg $(i,j) \in \mathcal{A}$ is given by $\delta_{(i,j)}$, and the processing time at a node $i \in \mathcal{S}$ is given by $p_i$. Each \distributioncenter, $i \in \mathcal{D}$, has a cutoff, $\phi_i$, that specifies the time by which a truck must have arrived for its contents to be processed same-day by the last-mile network.

The set of \paths through the network is given by $\ell \in \mathcal{L}$, where each \path is some sequence of legs $\ell = ((i, j), (j, k), \dots)$ with origin $o(\ell)$ and destination $d(\ell)$. Some \paths are \fc-\ds direct ($|\ell| = 1$) whereas others involve (potentially multiple) \scc's ($|\ell| > 1)$. The amount of volume carried by a \path is denoted by $vol(\ell)$ and the volume carried by a leg is $vol_{(i,j)}$ where $vol_{(i,j)} = \sum_{\ell \in \mathcal{L} : (i,j) \in \ell} vol(\ell)$. 

Let $\Pi(\ell)$ be the \textit{promise} of a \path $\ell \in \mathcal{L}$, representing the number of days that have elapsed between the \tdt at $o(\ell)$ and the arrival of the truck's contents at their ultimate destination (including the last mile portion). In general, the best possible promise for a \path is $\Pi(\ell) = 0$ (same day delivery), however, particularly for \paths involving long leg transition times, the promise can be significantly larger (e.g., a week or more).

All legs are expected, but not required, to have at least one \tdt placed on them. We discretize time into slots, where each slot represents a 15-minute time interval (yielding 96 slots over the course of a day). The set of slots that can have a \tdt placed is defined as $\mathcal{K} = \{1, \dots, 96\}$, such that $k \in \mathcal{K}$ for leg $(i,j)$ would mean that a truck departs from $i$ at slot $k$ en route to location $j$. Each leg can support some number of departure times called ``waves,'' $W_{(i,j)}$. In this study $W_{(i,j)} = 1$ for most legs but we consider multiple waves (specifically $W_{(i,j)} = 2)$ for \fc-\scc legs (i.e., those defined by $\{(i,j) \in \mathcal{A} : i \in \mathcal{F}, j \in \mathcal{S}\}$).\footnote{The set of waves for a leg is given by $\{1,\dots,W_{(i,j)}\}$ and we occasionally use the shorthand $[W_{(i,j)}]$ for compactness.}

\paragraph{Constraints} The problem has a variety of constraints that any \tdt placements must satisfy. (i) The first family of constraints are \textit{working shift} constraints where the set of slots that can support the inbound arrival of a truck (for nodes $i \in \mathcal{S} \cup \mathcal{D}$), or the outbound departure of a truck (for nodes $i \in \mathcal{F} \cup \mathcal{S}$), are specified as Boolean values. (ii) The problem also has both \textit{per-slot} and \textit{rolling capacity} constraints where, for each node in the network, limitations on the number of inbound and outbound trucks are expressed. (iii) Next, \textit{labor efficiency} constraints ensure that arriving volume is sufficiently distributed to mitigate potential backlog issues, and (iv) finally \textit{dispatch spacing} constraints ensure that there is some separation between \tdt values for multi-wave legs.

\paragraph{Objective function} We consider two objective functions in this work: the analytically expressible ``package speed'' of the network, and then a statistically-derived black-box objective that more accurately captures real-world dynamics. Given a set of assigned \tdt values on each leg we can calculate the package speed of a \path using the \tdt of the \path's first leg, $x_{(i,j)} \in \mathcal{K}$, and the promise of the \path:\footnote{Note that, for legs with multiple waves, the specific \tdt assigned to the \path must be specified in order to calculate package speed.}

\begin{equation}
f_{\text{package speed}}(\ell) = vol(\ell) \left(\frac{x_{(i,j)}}{|\mathcal{K}|} \cdot \omega(\Pi(\ell)) + (1-\frac{x_{(i,j)}}{|\mathcal{K}|}) \cdot \omega(\Pi(\ell)+1) \right)
\label{eqn:package-speed}
\end{equation}

\noindent where $\omega$ is a function that returns a weight for some promise value (higher weights are associated with lower promise). Intuitively, since the truck departs from the \fc at its \tdt, some portion of the volume for a \path is rewarded at promise $\Pi(\ell)$, while the remainder of the volume is rewarded at promise $\Pi(\ell)+1$. The overall package speed objective is then given by $F_{\text{package speed}} = \sum_{\ell \in \mathcal{L}} f_{\text{package speed}}(\ell)$.

While the package speed objective is convenient to express, it has a number of simplifying assumptions (e.g., a linear assignment of volume to each promise value) and may not reflect the true dynamics of the system. As such, for our second objective function, we assume access to a black-box evaluator that has been trained on a history of \tdts, \path promises, and the associated network speed. Then, for a new candidate assignment of \tdts, we can query the black-box evaluator to retrieve expected reward for different \path promises (0-day, 1-day, etc.) which can then be combined them into an aggregate network speed objective value. When viewed on a per-\path basis, and for some target promise day $\bar{\Pi}$, this black box evaluation is defined by:

\begin{equation} \label{eq:black-box}
f_{\text{black box}}(x, \ell, \bar{\Pi}) = 
     \begin{cases}
       0, &\quad\text{if $\Pi(\ell) > \bar{\Pi}$} \\
       f(x_{\ell[0]}), &\quad\text{if $\Pi(\ell) = \bar{\Pi}$} \\
       f(|\mathcal{K}|) &\quad\text{if $\Pi(\ell) < \bar{\Pi}$} \\
     \end{cases}
\end{equation}

\noindent The specific implementation details for this black-box function $f(x)$ are out of the scope of this paper. To calculate the network-level objective value \obj for a given solution $x$, we aggregate over all \paths and target promise days $\bar{\Pi} \in [0, 3]$, as follows:

\begin{equation} \label{eq:obj}
F_{\text{black box}}(x) = \sum_{\bar{\Pi},\ell}f_{\text{black box}}(x, \ell, \bar{\Pi})
\end{equation}

\begin{figure}[h]
  \centering
  \includegraphics[width=0.7\textwidth]{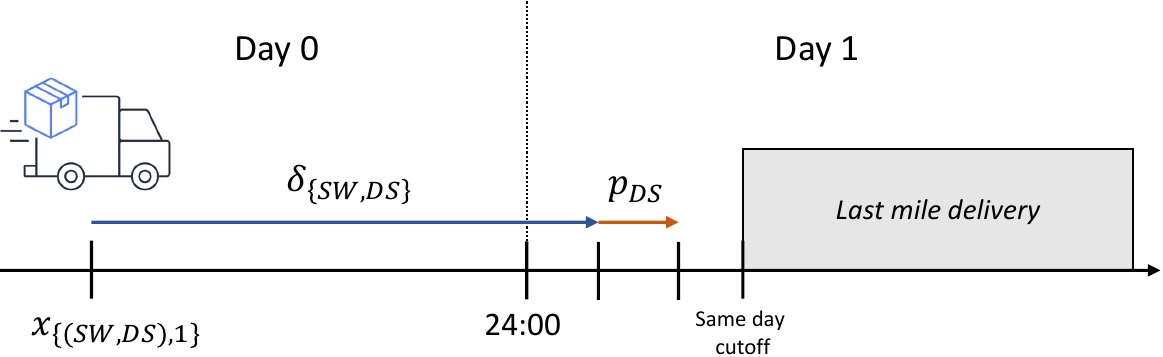}
  \caption{Example (direct) \fc-\ds connection. The \tdtime (\tdt) from the \fc node is $x_{(SW,DS),1}$ on Day 0 with transition time $\delta_{(SW,DS)}$ and processing time at the \ds of $p_{DS}$. If the truck arrival time ($x_{(SW,DS),1} + \delta_{(SW,DS)} + p_{DS}$) on the next day is prior to the same day cutoff time, the volume carried by the truck is delivered with next-day (i.e., 1-day) promise. If it is later than the same day cutoff time, it is delivered with 2-day promise. }
  \label{fig:next_day_delivery}
\end{figure}

\paragraph{Optimization task} The goal of the optimization problem is to produce a \tdt assignment to the legs that maximizes the objective function while satisfying problem constraints. Formally, a solution is defined as $X := \{(x_{(i,j),1},\dots, x_{(i,j),W_{(i,j)}}) : (i,j) \in \mathcal{A}\}$ where $x_{(i,j),w} \in \{\bot\} \cup \mathcal{K}$, with $\bot$ as a null value that indicates a \tdt is not placed on the leg. Thus, the goal of the optimization is to solve for $\max f(X) : \mathcal{C}$, where $f$ is a statistically derived black-box objective function and $\mathcal{C}$ is the set of problem constraints. In the following examples we detail a simple problem instance and how different \tdt assignments can cause \paths promise to vary.

\begin{example}[Direct \path]\label{example-1}
{Consider the direct \path defined by the leg $(SW2,DS3)$ in Figure \ref{fig:problem_schematic}}. Assume that $\delta_{(SW2, DS3)}$ is 12 hours, the cutoff time at DS3 is 8:00AM, and the processing time at DS3, $p_{DS3} = 0$. Any \tdt placed on the leg up to, and including, 8:00PM will result in a 1-day \path promise. A \tdt placed later than 8:00PM will result in a 2-day \path promise.
\end{example}

\begin{example}[Sort center \path]\label{example-2}
{Consider the \path $((SW1, SC1), (SC1, DS1))$  with leg transition times of 12 and 8 hours, respectively, and processing times of $p_{SC1} = 3$ and $p_{DS1} = 0$. DS1 has a cutoff of 8:00AM. A \tdt assignment of 9:00AM and 12:00AM on the legs, respectively, will yield a 1-day \path promise. However, a \tdt assignment of 9:00AM and 11:00PM on the legs, respectively, will yield a 2-day \path promise since the \tdt on the second leg occurs prior to the arrival of the truck from the first leg.} 
\end{example}

\section{Methods}
\label{sec:methods}

In this section we describe the core components of our CP and RKO solution approaches. For the CP approach we detail the decision variables, constraints, and objective function of our model. For the RKO approach we detail the anatomy of the algorithm, including the heuristic search logic and the decoder design. 

\subsection{Constraint programming (CP)} \label{sec:methods-cp}

The core decision variable in our CP model is $x_{(i,j),w} \in \mathcal{K}$, an integer variable that represents the \tdt assigned to leg $(i,j) \in \mathcal{A}$ for wave $w \in \{1,\dots,W_{(i,j)}\}$. We also introduce binary variables $y_{(i,j),k, w} \in \{0,1\}$ that take a value of $1$ if the \tdt for leg $(i,j) \in \mathcal{A}$ is $k$ for wave $w \in \{1,\dots,W_{(i,j)}\}$, and $0$ otherwise. 

The variables are related such that $y_{(i,j),k,w} = 1 \rightarrow x_{(i,j),w} = k$. Similarly, $y_{(i,j),k,w} = 0 \rightarrow x_{(i,j),w} \neq k$. We note that, for the purposes of this paper, we assume that $W_{(i,j)} = 1$ for all legs, except \fc-\scc legs (i.e., $(i,j) \in \mathcal{A}: i \in \mathcal{F}, j \in \mathcal{S}$) where $W_{(i,j)} = 2$. Further, if $W_{(i,j)} = 2$ and $x_{(i,j),1} = x_{(i,j),2}$ then the leg is considered to have one wave (since the \tdt variables are equal).

\paragraph{Working shifts constraints} Inbound and outbound shift constraints are embedded into the core variable domains. For a given $x_{(i,j),w}$ variable we remove the set of slots from $\mathcal{K}$ that do not permit an outbound departure from $i$ at that time. Similarly, we remove the set of slots from $\mathcal{K}$ that do not permit an arrival at $j$ at $k + \delta_{(i,j)}$. The resulting set of slots, $\mathcal{K}^{(i,j)}$, is then assigned as the domain, such that $x_{(i,j),w} \in \mathcal{K}^{(i,j)}$. Similarly, Boolean variables $y_{(i,j),k,w}$ are fixed to $0$ for all $k \in \mathcal{K} \setminus \mathcal{K}^{(i,j)}$, and untouched otherwise. 

\paragraph{Per-slot and rolling capacity constraints}
Inbound and outbound capacity constraints are expressed on a per-slot basis, as well as on a rolling slot basis. Let $c^{in}_j, c^{out}_j$ be the per-slot inbound and outbound capacities of node $j$, respectively. Per-slot constraints are as follows:
\begin{align}
\sum_{j \in \mathcal{S} \cup \mathcal{D}: (i,j) \in \mathcal{A}} \max(\{y_{(i,j),k,w} : w \in [W_{(i,j)}]\}) \leq c^{out}_i, &\forall k \in \mathcal{K}, i \in \mathcal{F} \cup \mathcal{S} \\ 
\sum_{i \in \mathcal{F} \cup \mathcal{S}: (i,j) \in \mathcal{A}} \max(\{y_{(i,j),k',w} : w \in [W_{(i,j)}]\}) \leq c^{in}_j, & \forall k \in \mathcal{K}, j \in \mathcal{S} \cup \mathcal{D} \label{ib-per-bucket}
\end{align}

\noindent where in Constraint~\eqref{ib-per-bucket},  $k' = (k - \delta_{(i,j)}) \bmod K$. For these constraints, the $\max$ operation ensures that legs that only have a single wave assigned do not double count capacity. 

Let $c^{in-rolling}_j, c^{out-rolling}_j$ be the inbound and outbound rolling capacities for node $j$, respectively. The rolling slot constraints include an additional sum over the rolling window, $\kappa$ and are defined as follows:

\begin{align}
\sum_{j \in \mathcal{S} \cup \mathcal{D}: (i,j) \in \mathcal{A}} \left(\sum_{k' \in \Omega^{OB}} \max(\{y_{(i,j),k,w} : w \in [W_{(i,j)}]\})\right) \leq c^{out-rolling}_i, & \forall k \in \mathcal{K}, i \in \mathcal{F} \cup \mathcal{S} \\ 
\sum_{i \in \mathcal{F} \cup \mathcal{S}: (i,j) \in \mathcal{A}} \left (\sum_{k' \in \Omega^{IB}} \max(\{y_{(i,j),k',w} : w \in [W_{(i,j)}]\}) \right) \leq c^{in-rolling}_j, &\forall k \in \mathcal{K}, j \in \mathcal{S} \cup \mathcal{D}
\end{align}

\noindent where $\Omega^{OB} = \{k, (k+1) \bmod K, \dots, (k+\kappa) \bmod K\}$ and $\Omega^{IB} = \{(k-\delta_{(i,j)}) \bmod K, (k - \delta_{(i,j)} + 1) \bmod K, \dots, (k - \delta_{(i,j)} + \kappa) \bmod K)\}$. 

\paragraph{Labor efficiency}

Our CP model includes labor efficiency constraints for each of the \distributioncenter nodes, $i \in \mathcal{D}$. The general modeling strategy is to: i) calculate the total volume on all of the legs inbound to a given \ds, ii) introduce a variable that represents the volume available for processing at the start of the shift, and then iii) introduce a series of variables that track the volume throughout the remainder of the shift, constrained appropriately.

Let $vol^{arriving}_j = \sum_{i \in \mathcal{V} : (i,j) \in \mathcal{A}} vol_{(i,j)}$ represent the total volume arriving at $j \in \mathcal{D}$ based on the volume of its inbound legs. Let variables $\alpha_{j, k} \in \{0,1, \dots, vol^{arriving}_j\}$ be used to track the volume remaining to be sorted at node $j$ in slot $k$. At the start of the shift for a particular \ds (e.g., $\text{start}_j$), the value of $\alpha_{j, \text{start}_j}$ can be determined by summing $y_{(i,j),k,w}$ for all legs inbound to $j$ and all $k$ between the cut-off for this DS, $\phi_j$, and the start of the shift at this DS, defined as follows:\footnote{In this work legs feeding into DS nodes only have a single wave.}

\begin{equation}
  \alpha_{j, \text{start}_j} = \sum_{i \in \mathcal{F} \cup \mathcal{S} : (i,j) \in \mathcal{A}}\sum_{k \in \Phi} vol_{(i,j)} \cdot y_{(i,j),k,1}, \forall j \in \mathcal{D}
\end{equation}

\noindent with $\Phi = \{(\phi_j + 1 - \delta_{(i,j)}) \bmod K, (\phi_j + 2 - \delta_{(i,j)}) \bmod K, \dots, \text{start}_j - 1 - \delta_{(i,j)} \bmod K\}$.

Next, we proceed to define the volume tracking variables for the remainder of the shift, $\alpha_{j,k} \forall k \in \text{shift} = [\text{start}_j + 1, \text{end}_j]$. Each tracking variable is defined based on its preceding tracking variable, as well as any new arrivals of volume as follows:\footnote{A small adjustment is made to the modeling for shifts that take place across the midnight boundary. Further, some adjustments are made in the model to ensure that negative backlogs are thresholded at zero.}

\begin{equation}
\alpha_{j, k} =  \alpha_{j, (k-1) \bmod K} - \lambda_j + \sum_{i \in \mathcal{F} \cup \mathcal{S} : (i,j) \in \mathcal{A}} vol_{(i,j)} \cdot y_{(i,j), k', 1}, \forall k \in \text{shift}_j, \forall j \in \mathcal{D}
\end{equation}

\noindent where $k' =(k-\delta_{(i,j)}) \bmod K$ and $\lambda_j = \frac{vol^{arriving}_j}{|\text{shift}_j|} $ is the daily volume processing rate for DS $j$ (adjusting for breaks where appropriate).

Finally, we constrain the final volume tracking variable to take on a value of $0$ (indicating that there is no additional backlog at the end of the shift) using:

\begin{equation}
    \alpha_{j, \text{end}_j} = 0, \forall j \in \mathcal{D}
\end{equation}

\noindent and also constrain the second last slot so it must be strictly greater than zero. 

\paragraph{Dispatch spacing} The dispatch spacing constraints, for legs involving multiple waves, ensure that there is some minimum spacing between each of the waves. Let $\epsilon_i \geq 0$ represent the dispatch spacing requirement for node $i \in \mathcal{F}$ (recalling that only \fc-\scc legs have multiple waves). We use the following constraints: 
\begin{equation}
\begin{split}
x_{(i,j),1} \neq x_{(i,j),2}  \rightarrow & |x_{(i,j),1} -x_{(i,j),2}| \geq \epsilon_i \wedge \\ & |(|K| + x_{(i,j),1}) -x_{(i,j),2}| \geq \epsilon_i \wedge 
\\ & |x_{(i,j),1} - (x_{(i,j),2} + |K|)| \geq \epsilon_i, 
\forall (i,j) \in \mathcal{F} \times \mathcal{S} 
\end{split}
\end{equation}

\noindent where the multiple terms enforce the required spacing while considering the wrapping of the planning horizon, |$\mathcal{K}$|. When $x_{(i,j),1} = x_{(i,j),2}$ the leg has one wave and spacing is not required.

\paragraph{Objective function}

Our CP approach is unable to include the black-box objective in a straightforward manner. As such, we seek to maximize the package speed equation, given by Eqn.~\eqref{eqn:package-speed}, in our CP model. In order to express this objective function, we need to encode the promise of each \path, $\Pi(\ell)$, in terms of the other variables in our model.\footnote{In our CP model, $\Pi(\ell)$ is represented as an integer variable.} As such, we iterate through each leg in the \path, $(i,j) \in \ell$, and introduce variables $\pi^{\ell}_{(i,j)}$ representing the additional promise incurred by the \tdt placement on that leg. The total promise for the path is then given by: 

\begin{equation}
    \Pi(\ell)  = \sum_{(i,j) \in \ell} \pi^{\ell}_{(i,j)}
\end{equation}

\noindent where $\pi^{\ell}_{(i,j)}$ is constrained based on the values of $x_{(i,j),w}$, $\delta_{(i,j)}$, $\phi_i$, and $p_i$. To ensure that only a single wave on each multi-wave leg is used to calculate promise for a given \ds, we introduce binary variables $x^{d}_{(i,j)}, \forall d \in \mathcal{D}$ and constrain them such that they take the value of one of the $x_{(i,j),w}$ variables. Then, $x^d_{(i,j)}$ is used in determining $\pi^{\ell}_{(i,j)}$ for the \ds $d$ associated with path $\ell$.

\paragraph{Implementation details}

We implement our CP approach using the CP-SAT solver from OR-Tools v9.11 (\cite{perron2019or}). The $\max$ expressions are implemented using an auxiliary variable and $\texttt{MaxEquality}$. For example, given an expression of the form $\max(vars) \leq B$, we introduce an auxiliary variable $aux = \texttt{MaxEquality}(vars)$ and then express $aux \leq B$. All floating point values (e.g., DS node volumes $\text{vol}_j$) are scaled so that the model can be expressed using all discrete values (integers), a requirement for working with CP-SAT. Logical implications (e.g., $\rightarrow$) are implemented using CP-SAT's \texttt{OnlyEnforceIf} functionality. Leg promise-added variables, $\pi^{\ell}_{(i,j)}$, used in formualting the objective function are expressed using a variety of auxiliary variables and constraints omitted here for brevity.  

\subsection{Random key optimizer (RKO)}\label{sec:methods-rko}

The RKO framework consists of two main components, comprising (i) problem-independent search logic, and (ii) problem-dependent decoding. 
See Algorithm \ref{alg:rko-search} for a formalized description. 
Any given \tdt assignment is encoded in the form of a random-key vector $\chi \in[0,1)^{\mathcal{N}}$, with $\mathcal{N} = \sum_{(i,j)\in\mathcal{A}}W_{(i,j)}$.
The problem-dependent RKO decoder $D(\chi)$ is a deterministic function that accepts such a random-key vector $\chi$ as input and outputs a \tdt solution candidate $X$ and its fitness value $f(X)$, i.e., $X, f(X)=D(\chi)$. 
As described in Algorithm \ref{alg:rko-search}, the fitness value $f(X)$ is used to guide the search process by iteratively updating the solution candidate $\chi$.  

\begin{algorithm}
\caption{RKO search process utilizing a basic annealing search engine}
\label{alg:rko-search}
\textbf{Input}: Problem metadata, hyperparameters $\theta$, initial solution $\chi_0$ \\
\textbf{Output}: Best solution found, best \tdt plan found, associated fitness
\begin{algorithmic}
\STATE Initialize \textsc{Decoder} $D(\chi) \mapsto x, f(x)$ \textcolor{gray}{\COMMENT{random key vector mapped to \tdt plan and fitness}} \\
\STATE Decode $D(\chi_0) \mapsto x_0, f(x_0)$ \\ 
\STATE Initialize $\chi_{best} \leftarrow \chi_0, x_{best} \leftarrow x_0$, $f_{best} \leftarrow f(x_0)$ \\
\STATE Initialize algorithm temperature $T_0$ \\
\FOR{$t=1$ through $t_{max}$ (inclusive)}
    \STATE $T_t \leftarrow \textsc{TemperatureUpdate}(t, T_{0})$  \\
    \STATE $\chi_t \leftarrow \textsc{Update}(\chi_{t-1})$ \textcolor{gray}{\COMMENT{propose new solution candidate}} \\
    \STATE $x_t, f_t \leftarrow D(\chi_t)$ \textcolor{gray}{\COMMENT{decode and evaluate fitness}}\\
    \STATE $\chi_t \leftarrow \textsc{Metropolis}(\chi_t, \chi_{t-1}, T_{t})$ \textcolor{gray}{\COMMENT{apply Metropolis acceptance}} \\
    Update current solution $\chi_{best}, x_{best}, f_{best}$
\ENDFOR
\RETURN $\chi_{best}, x_{best}, f_{best}$
\end{algorithmic}
\end{algorithm}

\paragraph{Random key encoding}\label{sec:methods-rko-encoding}

In our encoding the \tdt assignment for leg $(i,j)$ and wave $w$ is encoded by a single key $\chi_{(i,j),w}$, yielding the global random key representation $\chi = \bigcup_{(i,j)\in\mathcal{A}}\{\chi_{(i,j),w} : w \in \{1,\dots,W_{(i,j)]}\} \}$.
Given a fixed ordering of all legs $(i,j)$ in the network, we have a direct relationship between the $n$-th element of the random key vector $\chi$ and the leg and wave it represents, i.e., $n \leftrightarrow ((i,j),w)$, 
encoded in the map $M_{\mathrm{leg}}(n)$.

\paragraph{Decoder design}\label{sec:methods-decoder}

\begin{algorithm}
\caption{RKO decoder}
\label{alg:rko-decoder}
\textbf{Input}: Solution candidate $\chi$, \tdt bin map $M_{\mathrm{\tdt}}$, penalty weight $P$ \\
\textbf{Output}: Decoded \tdt plan $X$, associated solution fitness $f(X)$
\begin{algorithmic}
\FOR[\textcolor{gray}{for each random key in vector}]{$\chi_n \in \chi$}
    \STATE $x_n \leftarrow M_{\mathrm{\tdt}}(\chi_n)$ \textcolor{gray}{\COMMENT{map random key to associated \tdt value}}
\ENDFOR
\STATE $x_{paths} \leftarrow \textsc{Propagate}(X)$ \textcolor{gray}{\COMMENT{propagate \tdt  plan $X$ through network}}
\STATE $F_{\text{black box}} \leftarrow \textsc{Evaluate}(x_{paths})$ \textcolor{gray}{\COMMENT{black-box query for \obj}}
\STATE $N_{viol} \leftarrow \textsc{ViolationCheck}(x_{paths})$ \textcolor{gray}{\COMMENT{run constraint checking logic}}
\STATE $f(X) \leftarrow F_{\text{black box}} - P\cdot N_{\mathrm{viol}}$ \textcolor{gray}{\COMMENT{calculate solution fitness}}
\RETURN $X, f(X)$
\end{algorithmic}
\end{algorithm}

\begin{figure}[h]
  \centering
  \includegraphics[width=\textwidth]{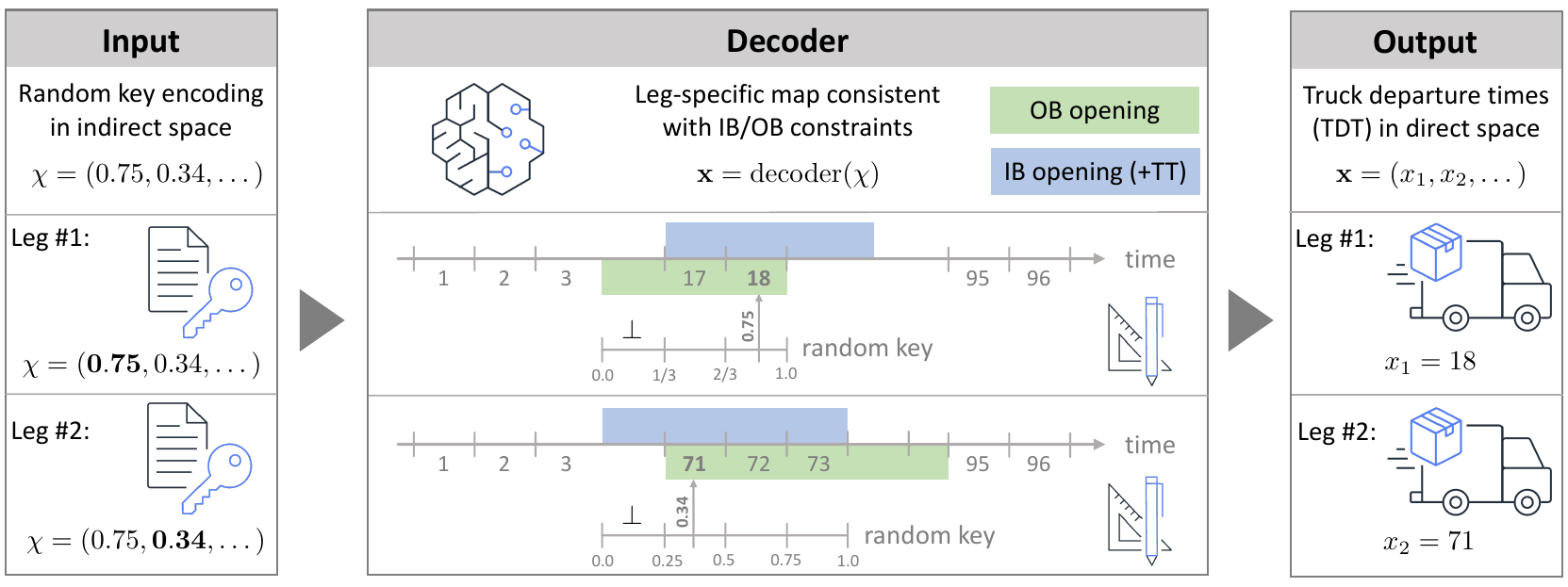}
  \caption{
  Schematic illustration of the RKO decoder design, with native handling of IB/OB constraints. 
  The input is given in terms of a random key $\chi$ with element $\chi_{n} \in [0,1)$ encoding the \tdt $x_{(i,j),w}$ for leg $n=(i,j)$. 
  The decoder maps $\chi_{n}$ to $x_{(i,j),w}$ while ensuring compatibility with the opening hours of the relevant inbound (IB) and outbound (OB) nodes. 
  To this end, the key $\chi_{n}$ is binned such that bins can only map to a \tdt assignment labeled as $\bot$ (where no \tdt is placed) or times compatible with OB opening hours and IB opening hours (adjusted accordingly by the corresponding travel time $\delta_{(i,j)}$), as indicated by the green and blue bars, respectively. 
  The decoder ensures mapping to IB/OB-feasible time domains only via leg-specific bins, as illustrated here for two different legs. 
  }
  \label{fig:rko_decoder_waves1}
\end{figure}

The decoder (Algorithm \ref{alg:rko-decoder}) accepts a random-key vector $\chi$ as input and outputs a complete \tdt assignment $X$ with fitness $f(X)$. 
As illustrated in Fig. \ref{fig:rko_decoder_waves1}, our decoder is designed such that compatibility with the opening hours of the relevant inbound (IB) and outbound (OB) nodes is ensured for every leg, adhering to the working shift constraint. 
This is done via a leg-specific mapping $\chi_{n} \rightarrow \{\bot\} \cup \mathcal{K}^{(i,j)}$, with $\mathcal{K}^{(i,j)}$ denoting the IB/OB-feasible time domains for leg $n=((i,j),w)$, which is contained in the map $M_{\tdt}$. 
In the following, we detail the subsequent decoding steps, labeled as (i) \textsc{Propagate}, (ii) \textsc{Evaluate}, and (iii) \textsc{ViolationCheck} in Algorithm \ref{alg:rko-decoder}.

\paragraph{Propagate}\label{sec:methods-rko-propagate}
In order to calculate the objective value \obj, and properly check for satisfaction of constraints, we must define a procedure for taking a given set of \tdt assignments and propagating them through the network. 
First, we pick a \path and, starting at the \fc $i$, we consider the set of \tdt values that could be assigned to that leg $\{x_{(i,j),w}\forall w \in W_{(i,j)}\}$.
We simulate the assignment of each \tdt option and calculate the expected promise ($\Pi_{w}(\ell)$) for each, for the \path.
To do this, we start from $x_{(i,j),w}$ and calculate the arrival time at the downstream node $j$.
We compare the arrival time, adjusted by node SLA $\mathrm{p}_j$, against the \tdt of the next leg $x_{(j,k)}$.
If we satisfy $x_{(i,j),w} + \delta_{(i,j)} + \mathrm{p}_j < x_{(j,k)}$ we assume the package will carry forward on the same day, otherwise the package must dwell and is assumed to carry forward on the next day.
The \ds nodes do not have an SLA; their arrival processing is dictated by the labor curve, as described in \ref{sec:problem}, and so to satisfy same-day delivery we must simply satisfy the cutoff $x_{(i,j)} + \delta_{(i,j)} < \phi_j$. 
The \tdt is then chosen as $\tdt_{(i,j),\ell} = \max_{(i,j),w}(\min_{(i,j),w}(\Pi_{w}(\ell))) \forall w \in W_{(i,j)}$. 
After propagating the \tdt solution across all \paths, we yield a set of delivery horizon promise days $\Pi(\ell)$ for all \paths, with one promise day per \path.

\paragraph{Evaluate}\label{sec:methods-rko-obj}
In order to query the demand covered by a \path, we take the \tdt associated with the origin node \fc and adjust for the \fc's order processing time to obtain an order cutoff time. 
We combine these cutoffs with the associated promise day $\Pi(\ell)$ for the \path we wish to query for.
We use these cutoffs and promise days $\Pi(\ell)$ to query a black box demand estimator, which yields an expected demand coverage for the \path (based on historical data) as described in Eqn.~\ref{eq:black-box}. 
However, as mentioned in Section \ref{sec:introduction}, because each \ds is serviced by multiple \fc's, and these \fc's have overlapping potential coverage, we must include all \paths that lead to a particular \ds, such that we can account for the shared demand coverage. 
Thus, we can calculate a reasonable demand coverage estimate on the downstream \ds level, but not at the \path level.
This demand coverage is summed across all \ds in the network, across all promise days of interest $\bar{\Pi}$, and summed to yield the top-level solution coverage, as shown in Eqn.~\ref{eq:obj}.

\paragraph{ViolationCheck}\label{sec:methods-rko-constraints}
Within the RKO framework, there will typically be a number of constraints that can be handled by design as part of the mapping procedure from random key to problem domain, which we refer to as “native” or “implicit” handling, and there will be a number of constraints that cannot be satisfied natively and, as such, must be addressed via explicit penalty terms.
In our case, we are able to (implicitly) satisfy the IB and OB shift constraints (constraint 1) by first filtering the possible decoding time buckets by the intersection of the leg origin's OB shift times and the leg destination’s IB shift times, yielding a reduced set of candidate \tdt bins, which must satisfy the downstream IB shift constraint. 

The IB and OB capacity constraints (constraint 2), both on the per-time-bucket and rolling window basis, are simple to calculate: we count the number of unique time bins within a rolling time window of a certain size. 
For the per-time-bucket check, this window size is one; for the rolling hour basis, the window size is four. 
Similarly, for dispatch spacing (constraint 4), we can calculate the difference between unique \tdts at a node, and compare that against the specified threshold value.

For labor efficiency (constraint 3), we start by calculating all the arrivals from within the working shift, and anything arriving outside of that window is considered starting backlog for the day (anything arriving after shift end must be assigned to the next promise day). 
All backlog and arriving volume require one time bucket to be prepared to be processed by available labor. 
For every active shift time bin, we subtract the labor processing volume from the current backlog to yield the future backlog. 
If the future backlog is not zero by the end of the shift window, we count and report a constraint violation.

\paragraph{Fitness} 
After the RKO decoder has decoded and propagated a given solution candidate $\chi \mapsto X$, it must calculate the associated fitness score $f(X)$.
The associated \obj score and $N_{viol}$ count are calculated for the solution $X$, as described in \ref{sec:methods-rko-obj} and \ref{sec:methods-rko-constraints}.
To generate a signal for the searcher to learn with, instead of simply marking any solutions where $N_{viol} > 0$ with some (infinitely) large penalty for infeasibility, we choose some more moderate penalty factor $P$.
This way, the search algorithm receives stepwise feedback as it reduces the number of violations in its proposed solutions. 
Thus, the fitness score is calculated as:

\begin{equation}
\label{eq:rko-fitness}
f(\chi) = F_{\text{black box}}(\chi) - P \cdot N_{\mathrm{viol}}(\chi)
\end{equation}

\paragraph{Law of physics solution}\label{sec:methods-lop}
An important concept in this problem is that of the law of physics (LOP) solution. 
By this we mean the set of \tdts along a \path that would minimize the time it takes for volume to travel that \path.
This LOP is a sort of theoretical "best" value for the \tdt, in that it establishes the absolute latest time you can place a \tdt and still make some target delivery horizon cutoff. 
Of course, while this is true on the per-\path basis, it may not be the case when consolidating into, for instance, a single \tdt value for a leg, to service all \paths along it.

We calculate the network LOP times by starting at a \ds and assuming the volume must arrive by the cutoff time $\phi$.
From there, we step back along each leg $(i,j)$ feeding to that \ds, and account for the travel time $\delta_{(i,j)}$ of each. 
For each \scc in the chain, we must account for the associated SLA $p_i$ time, but we repeat the process of stepping back through the network and adjusting for $\delta$ and $p$.
Once we hit the \fc, we are done and have our full set of LOP \tdts for all legs considered.
Where there are legs with multiple LOP \tdts, we pick the option(s) that would minimize the promise day $\Pi$, and if there are still multiple, we pick the latest option.
While this may not be optimal logic, it does ensure we are maximizing the volume that would arrive at the minimal promise day allowed by the joint set of \paths. 
This LOP solution can serve as a useful starting point for solution searches, even though the full solution can be (and is, in this case) infeasible, per the problem constraints.

\paragraph{RKO warm start} 
We can also encode a given solution $X$ into a random key vector $\chi$, by pseudo-inverting the decoder logic. 
Specifically, for a given \tdt solution vector $X$, we step through time assignments $x_{(i,j),w}$, identify the appropriate bin edges for the appropriate leg from $M_{\mathrm{\tdt}}(\cdot)$, and map the time to a random key $\chi_{k}$ that is within the range of the corresponding bin edges.
This enables seeding (warm-starting) the RKO solver with a solution from a different solver or process, potentially expediting the search or yielding improved results.

\section{Experiments}
\label{sec:experiments}

Here we describe our experimental setup and key results.
We start by introducing a baseline solver, followed by (hybrid) extensions of our solvers described in Section \ref{sec:methods}.
We then introduce the full problem instance, and describe our figures of merit.
Finally, we present the results of all solvers (baseline, RKO, CP) on the problem instance for ``single wave'' and ``multi wave'' problem variants.

\paragraph{Baseline approach}
\label{sec:methods-incumbent}

The baseline is a simple yet powerful approach for \tdt optimization, consisting of a greedy algorithm and subsequent local search procedure. 
The greedy algorithm proceeds by leg type. 
It first assigns \tdts on \fc-\ds and \scc-\ds legs before moving to indirect legs, i.e., legs that end in a \scc. 
For legs not starting at an \fc, the greedy algorithm proxies \obj since the \fc \tdt might not be available yet. 
For each leg type, the greedy algorithm iteratively finds the best feasible \tdt candidate for every leg based on \obj and places a \tdt on the leg that, given all previously placed \tdts, provides the largest incremental \obj benefit. 
Upon completion of the greedy algorithm, a local search procedure is launched for further solution refinement. 
This routine iteratively enumerates and applies \tdt changes to the initial solution. 
Candidates include single or multiple concurrent \tdt changes and are constructed using domain-specific knowledge of \tdt optimization to reduce the search space. 
The local search heuristic to select among such candidates is a classical hill-climbing algorithm.

\paragraph{Hybrid solver}
\label{sec:methods-hybrid}

We also test a hybrid approach where RKO and CP solutions are used to seed the baseline's local search routine. 
We label experiments that include this local search hybridization as ``hybrid'' (with ``baseline'' indicating standalone experiments).

\subsection{Problem instances and figures of merit}
\label{sec:experiments-instances}

\paragraph{Base instance} We take a single real-world, EU-scale middle mile network as our core problem instance. 
The network spans ten countries in the EU. 
The problem instance comprises $90$ \fc, $34$ \scc, and $242$ \ds nodes, with $2187$ legs and $11830$ \paths ($1010$ direct and $10820$ indirect \paths, respectively). 
All results reported in Section \ref{sec:results-problem} are based on this problem instance. 

\paragraph{Figures of merit} Solver results are characterized in terms of a set of relevant KPIs, specifically the 0D, 1D, 2D, and 3D demand coverage, and the number of \tdts placed. 
These KPIs are a measure of how much total network demand has been covered by the solution up to each target promise. 
This is equivalent to the value of \obj for each target promise day $\bar{\Pi}$, normalized by a global constant representing the maximum possible coverage of the network, such that the values are bounded in $[0, 1]$.
The KPI values are monotonically increasing, as the demand covered by a solution on target promise day $\bar{\Pi}$ also necessarily includes the demand covered on promise day $\bar{\Pi} - 1$, assuming $\bar{\Pi} > 0$.
We also calculate the difference between each solver’s KPI and the baseline, which serves as the benchmark result. 
We report the differences in terms of basis points (bps), where 1 bps is equivalent to a demand change of 0.0001, and we round to the nearest 10 bps (0.001 change in demand). At the scale of large e-commerce retailers, small values of bps improvement are worth multiple million dollars of speed benefit to customers.

Particular emphasis is placed on maximizing 1D demand coverage, following our discussion in Section \ref{sec:introduction}, so we use this as our primary decision criteria when considering the best solution strategy. 
We break out results by problem context, ``single wave'' with a maximum number of waves of one, and ``multi wave'' with a maximum number of waves of two on \fc-\scc legs. 

\subsection{Results}
\label{sec:results-problem}

\begin{table}[h]
\centering
\begin{tabular}{ c  c  c  c  c  c  c }
\toprule
 & Baseline & RKO  & RKO (hybrid) & CP  & CP (hybrid) \\ [0.5ex]
\midrule
\# \tdt placed & 2173 & 2182 & 2182 & 2187 & 2187 \\
\midrule
0D KPI & 0.005 & 0.007 & \bf{0.008} & 0.006 & \bf{0.008} \\
1D KPI & 0.522 & 0.515 & 0.526 & 0.522 & \bf{0.529} \\
2D KPI & \bf{0.712} & \bf{0.712} & 0.711 & 0.702 & 0.707 \\
3D KPI & 0.762 & \bf{0.765} & \bf{0.765} & 0.757 & 0.759 \\
\midrule
0D KPI $\pm$ (bps)	& 0 & 20  & 30  & 10   & 30 \\
1D KPI $\pm$ (bps) & 0 & -70 & 40  & 0    & 70 \\
2D KPI $\pm$ (bps) & 0 & 0   & -10 & -100 & -50 \\
3D KPI $\pm$ (bps) & 0 & 30  & 30  & -50  & -30 \\
\bottomrule
\end{tabular}
\caption{Results for solvers run on ``single wave'' variant (e.g., $W_{(i,j)} = 1, \forall (i,j) \in \mathcal{A}$). For reference, we include the number of \tdts placed by each solver. Our key metrics are the KPIs, which describe how much expected demand was covered according to each target promise, as a fraction of an estimated total network demand. We also present diffs between each solver's KPI and the baseline. The hybrid approaches combining RKO or CP with local search are found to perform the best with respect to the 0D and 1D \obj coverage KPI, though that seems to come at the cost of 2D coverage. 
Both the RKO and CP solvers are found to place \tdt assignments on more legs than the baseline solver, with the CP solver placing \tdts on all $2187$ legs, by design.}
\label{table:results-v0-kpis}
\end{table}

In Tab.~\ref{table:results-v0-kpis}, we present the best result from each approach on ``single wave,'' including the baseline solver as reference, the RKO standalone, RKO with local search hybrid, CP standalone, and CP with local search hybrid. 
The RKO standalone result is sourced from the best hybrid run, meaning it may not be the best RKO standalone solution, but it is the best RKO hybrid solution. 
The same applies to the CP results. 

We make the following observations.
First, both RKO and CP, even without local search, place more \tdts than the baseline, with the CP solver placing them on all $2187$ legs, by design. 
The requirement that CP place all \tdts available in the ``single wave'' variant was an assumption included to simplify the model, as described in Section~\ref{sec:methods-cp}.
It is generally expected that placing more \tdts will yield more demand coverage, but due to the non-linear objective (see Section~\ref{sec:introduction}), there is a limit to this effect. 
Further, due to the complex constraints, it can actually be beneficial to avoid placing some \tdts in order to yield more objective value out of the network, as these results suggest.
Second, we see that both the RKO and CP baseline, i.e., without local search, struggle to compete with the baseline solution on 1D coverage. 
However, we also see that both hybrid approaches outperform the baseline with respect to 0D and 1D coverage, with $+40$ bps and $+70$ bps improvements in next-day coverage for RKO-hybrid and CP-hybrid, respectively. 
It is worth noting that small bps improvement ($+40$ bps) can yield multiple millions USD in improved customer order volumes to large e-commerce retailers.
Interestingly, this seems to come at the expense of 2D, which may or may not be acceptable to the business.

\begin{table}[h]
\centering
\begin{tabular}{ c c c c c c c }
\toprule
 & Baseline & RKO & RKO (hybrid) & CP & CP (hybrid) \\ [0.5ex]
\midrule
\# \tdt placed & 2690 & 2763 & 2763 & 2721 & 2721 \\
\midrule					
0D KPI & 0.005 & 0.005 & \bf{0.007} & 0.005 & \bf{0.007} \\
1D KPI & 0.530 & 0.527 & \bf{0.535} & 0.523 & 0.532 \\
2D KPI & \bf{0.721} & 0.720 & 0.719 & 0.713 & 0.714 \\
3D KPI & 0.767 & \bf{0.768} & \bf{0.768} & 0.762 & 0.764 \\
\midrule					
0D KPI $\pm$ bps & 0 & 0   & 20  & 0   & 20 \\
1D KPI $\pm$ bps & 0 & -30 & 50  & -70 & 20 \\
2D KPI $\pm$ bps & 0 & -10 & -20 & -80 & -60 \\
3D KPI $\pm$ bps & 0 & 10  & 10  & -50 & -40 \\
\bottomrule
\end{tabular}
\caption{Results for solvers run on ``multi wave'' variant, i.e., \fc-\scc legs have two waves, all other legs have one. The hybrid approaches combining RKO or CP and local search perform the best with respect to the 0D and 1D KPI, although CP is found to be less competitive compared to the ``single wave'' results above. Again both RKO and CP solvers place more \tdts on more legs than the baseline.}
\label{table:results-v1-kpis}
\end{table}

In Tab.~\ref{table:results-v1-kpis}, we present the best result from each approach on the ``multi wave'' problem, including the same solver set as in Tab.~\ref{table:results-v0-kpis}, and the same relationship between best hybrid result and the associated baseline. 
We see similar trends as with the ``single wave'' results. Without local search, both RKO and CP solvers struggle to compete with the baseline approach, but when incorporating that local search phase, the hybrid solvers outperform the baseline in 0D and 1D coverage, again coming at the expense of 2D coverage. 
Concretely, the hybrid RKO approach yields $+50$ bps, and hybrid CP yields $+20$ bps over the baseline for 1D coverage. 
We note that the RKO results are derived from a warm start run, and so in aggregate it required multiple days of runtime, whereas the CP result was from a single run, with a budgeted runtime of a few hours.
Thus it appears that RKO is more competitive than CP in this setting, although the difference is small in 0D and 1D, but more pronounced in 2D and 3D.

\section{Conclusions and Outlook}
\label{sec:conclusion}

In this work we have designed and implemented algorithms for solving the middle-mile next-day delivery optimization problem at business-relevant scales, adopting two complementary approaches: the exact CP framework, and the heuristic RKO framework. 
We have studied the problem in two settings, one where we solve the problem assuming one \tdt per leg (``single wave''), and one where we allow for up to two \tdts on \fc-\scc legs (``multi wave''), which also required additional dispatch spacing constraints. In particular, we have shown that hybrid solvers integrating CP or RKO with custom local search can outperform the baseline solution in expected next-day delivery coverage. 

In future work it will be interesting to consider expanded problem definitions (e.g., in the form of additional features and constraints) and to explore alternative algorithmic strategies. 
For example, for RKO one could investigate alternative decoder designs or additional RKO optimization engines~\citep{chaves:24rko}.
Further, one could use path-relinking techniques for potential solution refinements. On the CP side, one could investigate decomposition methods such as CP-based column generation~\citep{junker1999framework}, custom inference techniques (e.g., problem-specific global constraints), or a tighter integration of the black-box objective into the search.  

\section*{Acknowledgements}

The authors would like to thank Helmut Katzgraber for his support in this work through his rigorous feedback and various insightful discussions. The authors would also like to thank Michal Penarski for his assistance in managing this project.

\bibliographystyle{elsarticle-num}
\bibliography{main}

\section{Appendix}

\paragraph{Synergies}
\label{sec:methods-synergies}

Pursuing parallel modeling branches allows us to explore potential interplay between those branches. 
For instance, one interesting feature of the CP approach is that it can output a set of pruned domains for each variable according to the constraints, effectively trimming the search space. 
This is an important part of the CP internal workings that RKO does not attempt. 
However, because we explicitly define the domain of each leg for the sake of time binning within the RKO decoder design, we can incorporate the CP pruned domains to help make the RKO search more efficient. 
Specifically, we have implemented this hybrid strategy in the context of IB/OB shifts, with the RKO decoder only decoding keys into domains as pruned by CP. 

Yet another, straightforward interplay is to take the output of one algorithm and use it to seed the starting point of another. 
This is typically referred to as “warm starting” a solver. 
We adopt this strategy as a basis for the hybrid solver, as described in Section \ref{sec:methods-hybrid}, but can also do this between the CP and RKO solvers. 
Additionally, we can warm restart the RKO solver by seeding it with an output from a different RKO run. 
This can be beneficial, as it effectively restarts the solver dynamics, yet at a higher quality solution while effectively resetting the temperature to high values that favor exploration.

\end{document}